\documentclass[a4, 12pt]{amsart}
\usepackage{amssymb}
\usepackage{amstext}
\usepackage{amsmath}
\usepackage{amscd}
\usepackage{latexsym}
\usepackage{amsfonts}
\usepackage[all]{xy}

\theoremstyle{plain}
\newtheorem{Theorem}{Theorem}[section]

\newtheorem{thm}[Theorem]{Theorem}

\newtheorem{prop}[Theorem]{Proposition}

\newtheorem{lem}[Theorem]{Lemma}

\newtheorem{cor}[Theorem]{Corollary}

\newtheorem{conj}[Theorem]{Conjecture}

\newtheorem*{discussion*}{Discussion}

\newtheorem*{thm*}{Theorem}
\newtheorem*{cor*}{Corollary}

\newtheorem*{claim*}{Claim}

\theoremstyle{definition}

\newtheorem{example}[Theorem]{Example}

\newtheorem{remark}[Theorem]{Remark}

\newtheorem{question}[Theorem]{Question}

\theoremstyle{remark}

\newcommand{\Proof}{\begin{proof}}
\newcommand{\Qed}{\end{proof}}

\newcommand{\rma}{\mathrm{a}}

\newcommand{\rme}{\mathrm{e}}

\newcommand{\rmr}{\mathrm{r}}

\newcommand{\rmF}{\mathrm{F}}
\newcommand{\rmG}{\mathrm{G}}

\newcommand{\calR}{\mathcal{R}}

\newcommand{\fka}{\mathfrak{a}}

\newcommand{\fkm}{\mathfrak{m}}

\newcommand{\fkq}{\mathfrak{q}}

\def\Ht{\operatorname{ht}}

\def\Spec{\operatorname{Spec}}

\def\GCD{\operatorname{GCD}}

\begin{document}

\setlength{\baselineskip}{20pt}

\title{Quasi-socle ideals in local rings 
with Gorenstein tangent cones}
\author{Shiro Goto}
\address{Department of Mathematics, School of Science and Technology, Meiji University, 1-1-1 Higashimita, Tama-ku, Kawasaki 214-8571, Japan}
\email{goto@math.meiji.ac.jp}
\author{Satoru Kimura}
\address{Department of Mathematics, School of Science and Technology, Meiji University, 1-1-1 Higashimita, Tama-ku, Kawasaki 214-8571, Japan}
\email{skimura@math.meiji.ac.jp}
\author{Naoyuki Matsuoka}
\address{Department of Mathematics, School of Science and Technology, Meiji University, 1-1-1 Higashimita, Tama-ku, Kawasaki 214-8571, Japan}
\email{matsuoka@math.meiji.ac.jp}
\author{Tran Thi Phuong}
\address{Department of Information Technology and Applied Mathematics,
Ton Duc Thang University, 98 Ngo Tat To Street, Ward 19, Binh Thanh District,
Ho Chi Minh City, Vietnam}
\email{sugarphuong@gmail.com}
\thanks{
\endgraf
}
\maketitle
\begin{abstract}
Quasi-socle ideals, that is the ideals $I$ of the form $I= Q : \fkm^q$ in a Noetherian local ring $(A, \fkm)$ with the Gorenstein tangent cone ${\rmG}(\fkm) = \bigoplus_{n \geq 0}{\fkm}^n/{\fkm}^{n+1}$  are explored, where $q \geq 1$ is an integer and $Q$ is a parameter ideal of $A$ generated by monomials of a system $x_1, x_2, \cdots , x_d$ of elements in $A$ such that $(x_1, x_2, \cdots , x_d)$ is a reduction of $\fkm$. 
The questions of when $I$ is integral over $Q$ and of when the graded rings $\rmG(I) = \bigoplus_{n \geq 0}I^n/I^{n+1}$ and $\rmF(I) = \bigoplus_{n \ge 0}I^n/\fkm I^n$ are Cohen-Macaulay are answered. Criteria for $\rmG (I)$ and $\calR (I) = \bigoplus_{n \geq 0}I^n$ to be Gorenstein rings are given. 
\end{abstract}
\section{Introduction}

This paper aims at a study of quasi-socle ideals in a local ring with the Gorenstein tangent cone. Our purpose is to answer Question \ref{1.1} below, of when the graded rings associated to the ideals are Cohen-Macaulay and/or Gorenstein rings, estimating their reduction numbers with respect to minimal reductions.

Let $A$ be a Noetherian local ring with the maximal ideal $\fkm$ and $d=\dim A > 0$.
Let $Q=(x_1, x_2, \cdots ,x_d)$ be a parameter ideal in $A$ and let $q \geq 1$ be an integer. We put $I=Q:\fkm^q$ and refer to those ideals as quasi-socle ideals in $A$. Then one can ask the following questions, which are the main subject of the researches \cite{GMT, GKM} and the present research as well.

\begin{question}\label{1.1}
\item[$(1)$] Find the conditions under which $I \subseteq \overline{Q}$, where $\overline{Q}$ stands for the integral closure of $Q$.
\item[$(2)$] When $I \subseteq \overline{Q}$, estimate or describe the reduction number 
$$\rmr_Q(I) = \min~\{0 \leq n \in \Bbb Z \mid I^{n+1} = QI^n \}$$
of $I$ with respect to $Q$ in terms of some invariants of $Q$ or $A$.
\item[$(3)$] Clarify what kind of ring-theoretic properties of the graded rings associated to the ideal $I$
$$\calR (I) = \bigoplus_{n \geq 0}I^n, ~~ \rmG (I) = \bigoplus_{n \geq 0}I^n/I^{n+1}, ~~ \operatorname{and}~~\rmF (I) = \bigoplus_{n \geq 0}I^n/\fkm I^n$$
enjoy.
\end{question}

In this paper we shall focus our attention on a certain special kind of quasi-socle ideals. We now assume that the tangent cone, that is the associated graded ring $$\rmG(\fkm) = \bigoplus_{n \ge 0 } \fkm^n / \fkm^{n+1}$$ of $\fkm$, is a Gorenstein ring and that the maximal ideal $\fkm$ contains a system $x_1, x_2, \cdots , x_d$ of elements such that the ideal $(x_1, x_2, \cdots , x_d)$ is a reduction of $\fkm$ (the latter condition is always satisfied if the field $A/\fkm$ is infinite; see \cite{NR} for the existence of reductions of $\fkm$ generated by $d$ elements). 
Let $a_1, a_2, \cdots , a_d,$ and $q$ be positive integers and we put
$$
Q=(x_1^{a_1}, x_2^{a_2}, \cdots , x_d^{a_d}) \ \ \textup{and} \ \ I=Q:\fkm^q.
$$
Let $\overline{A} = A/Q$, $\overline{\fkm} = \fkm /Q$, and $\overline{I} = I/Q$. Let $\rho = \max~\{n \in \Bbb Z \mid \overline{\fkm}^n \ne (0)\},$
that is index of nilpotency of the ideal $\overline{\fkm}$, and put $$\ell = \rho + 1 - q.$$
We then have the following, which are the answers to Question \ref{1.1} in our specific setting.

\begin{thm}\label{1.2}
The following three conditions are equivalent to each other.
\begin{enumerate}
\item[$(1)$] $I \subseteq \overline{Q}$. 
\item[$(2)$] $\fkm^qI = \fkm^qQ$.
\item[$(3)$] $\ell \geq a_i$ for all $1 \leq i \leq d$.
\end{enumerate}
When this is the case, the following assertions hold true.
\begin{enumerate}
\item[{\rm (i)}] $\rmr_Q(I) = \lceil \frac{q}{\ell} \rceil := \min~\{ n \in \Bbb Z  \mid  \frac{q}{\ell} \leq n \}$.
\item[{\rm (ii)}] The graded rings  $\rmG (I)$ and $\rmF (I)$ are Cohen-Macaulay.
\end{enumerate}
\end{thm}

\begin{thm}\label{1.2'}
Suppose that $\ell \geq a_i$ for all $1 \leq i \leq d$. Then we have the following. 
\begin{enumerate}
\item[{\rm (i)}] $\rmG (I)$ is a Gorenstein ring if and only if $\ell \mid q$.
\item[{\rm (ii)}] $\calR (I)$ is a Gorenstein ring if and only if $q  = (d-2)\ell$. 
\end{enumerate}
\end{thm}

Our setting naturally contains the case where $A$ is a regular local ring with $x_1, x_2, \cdots , x_d$ a regular system of parameters, the case where $A$ is an abstract hypersurface with the infinite residue class field, and the case where $A = R_M$ is the localization of the homogeneous Gorenstein ring $R= k[R_1]$ over an infinite field $k = R_0$ at the irrelevant maximal ideal $M = R_+$. In Section 3 we will explore a few examples, including these three cases,  in order to see how Theorems \ref{1.2} and \ref{1.2'} work for the analysis of concrete examples. The proofs of Theorems \ref{1.2} and \ref{1.2'} themselves shall be given in Section 2. However, before entering details, let us here explain the reason why we are interested in quasi-socle ideals.

The study of $\it{socle}$ ideals dates back to a research of L. Burch \cite{B}, where she explored certain socle ideals of finite projective dimension and gave a characterization of regular local rings (cf. \cite[Theorem 1.1]{GH}). More recently,
A. Corso and C. Polini \cite{CP1, CP2} studied, with interaction to the linkage theory of ideals, the socle ideals $ I = Q : \fkm$ in a Cohen-Macaulay local ring $(A, \fkm)$ and showed that $I^2 = QI$, once $A$ is {\it not} a regular local ring. Consequently, the rings $\rmG (I)$ and $\rmF (I)$ are Cohen-Macaulay and so is the ring $\calR (I)$ if $\dim A \geq 2$. The first author and H. Sakurai \cite{GSa1, GSa2, GSa3} explored also the case where the base ring is not necessarily Cohen-Macaulay but Buchsbaum, and showed that the equality $I^2 = QI$ (here $I = Q : \fkm$) holds true for numerous parameter ideals $Q$ in a given Buchsbaum local ring $(A, \fkm)$, whence $\rmG (I)$ is a Buchsbaum ring, provided that $\dim A \geq 2$ or that  $\dim A = 1$ but the multiplicity $\rme_\fkm^0(A)$ of $A$ with respect to $\fkm$ is not less than 2.

A more important thing is, however, the following. If $J$ is an equimultiple Cohen-Macaulay ideal of reduction number one in a Cohen-Macaulay local ring, the associated graded ring $\rmG(J) = \bigoplus_{n \ge 0}J^n/J^{n+1}$ of $J$ is a Cohen-Macaulay ring and, so is the Rees algebra $\calR(J)  = \bigoplus_{n \ge 0} J^n$ of $J$, provided $\Ht_A J \geq 2$. One knows the number and degrees of defining equations of $\calR(J)$ also, which makes the process of desingularization of $\Spec A$ along the subscheme $\mathrm{V}(J)$ fairly explicit to understand. 
This observation motivated the ingenious research of C. Polini and B. Ulrich \cite{PU}, where they posed, among many important results, the following conjecture.

\begin{conj}[\cite{PU}]\label{1.3}
\textit{Let $(A, \fkm)$ be a Cohen-Macaulay local ring with $\dim A \geq 2$. Assume that $\dim A \geq 3$ when  $A$ is regular. Let $q \geq 2$ be an integer and  let $Q$ be  a parameter ideal in $A$ such that $Q \subseteq \fkm^q$. Then
 $$Q:\fkm^q \subseteq \fkm^q.$$ }
\end{conj}

This conjecture was settled by H.-J. Wang \cite{Wan}, whose theorem says:

\begin{thm}[\cite{Wan}]\label{1.4}
Let $(A, \fkm)$ be a Cohen-Macaulay local ring with $d = \dim A \ge 2$. Let $q \ge 1$ be an integer and $Q$  a parameter ideal in $A$.
Assume that $Q \subseteq \fkm^q$ and put $I=Q:\fkm^q$.
Then $$I\subseteq \fkm^q, \ \ \fkm^q I = \fkm^q Q,\ \ ~~\operatorname{and}~~\ \ I^2=QI,$$ provided that $A$ is not regular if $d \ge 2$ and that $q \ge 2$ if $d \ge 3$.
\end{thm}

The recent research \cite{GMT} of the first and the third authors jointly with R. Takahashi  reports a different approach to the Polini-Ulrich Conjecture \ref{1.3} and has proven the following.

\begin{thm}[\cite{GMT}]\label{1.5}
Let $(A, \fkm)$ be a Gorenstein local ring with $d=\dim A>0$ and $\rme^0_\fkm(A) \ge 3$, where $\rme^0_\fkm(A)$ denotes the multiplicity of $A$ with respect to $\fkm$.
Let $Q$ be a parameter ideal in $A$ and put $I=Q:\fkm^2$.
Then $\fkm^2 I = \fkm^2 Q$, $I^3=QI^2$, and $\rmG(I)= \bigoplus_{n \geq 0}I^n/I^{n+1}$ is a Cohen-Macaulay ring, so that $\calR (I)=\bigoplus_{n \geq 0}I^n$ is also a Cohen-Macaulay ring, provided $d \ge 3$.
\end{thm}

The researches \cite{Wan} and \cite{GMT} are performed  independently and their methods of proof are totally different from each other's. The technique of \cite{GMT} can not go beyond the restrictions that $A$ is a Gorenstein ring, $q=2$, and $\rme_\fkm^0(A) \ge 3$. However, despite these restrictions, the result \cite[Theorem 1.1]{GMT} holds true even in the case where $\dim A=1$, while Wang's result says nothing about the case where $\dim A=1$. As is suggested in \cite{GMT}, the one-dimensional case is substantially different from higher-dimensional cases and more complicated to control. This observation has led the first three authors to the research \cite{GKM}, where they have explored quasi-socle ideals in Gorenstein numerical semigroup rings over fields as the first step towards further investigations of the Polini-Ulrich Conjecture \ref{1.3} of arbitrary dimension. The present research is, more or less, a continuation of \cite{Wan, GMT, GKM},  and our Theorems \ref{1.2} and \ref{1.2'} might have some significance towards the second step, providing some insight with Question \ref{1.1}, that are not presented by \cite{Wan, GMT, GKM}.

\section{Proof of Theorems \ref{1.2} and \ref{1.2'}}
The purpose of this section is to prove Theorems \ref{1.2} and \ref{1.2'}. First of all, let us restate our setting, which we shall maintain throughout this paper.

Let $A$ be a Noetherian local ring with the maximal ideal $\fkm$ and $d=\dim A > 0$.
We assume that the associated graded ring $$\rmG(\fkm) = \bigoplus_{n \ge 0 } \fkm^n / \fkm^{n+1}$$ of $\fkm$ is Gorenstein and that the maximal ideal $\fkm$ contains a system $x_1, x_2, \cdots , x_d$ of elements which generates  a reduction  of $\fkm$. Then, since $\rmG (\fkm)$ is a Gorenstein ring, the base local ring  $A$ is also Gorenstein and the initial forms $\{X_i\}_{1 \le i \le d}$ of $\{x_i\}_{1 \leq i \leq d}$ with respect to $\fkm$ constitute a regular sequence in $\rmG(\fkm)$, so that we get a canonical isomorphism 
$$\rmG(\fkm/(x_1, x_2, \cdots,x_d)) \cong \rmG(\fkm)/(X_1, X_2, \cdots, X_d)$$
of graded $A$-algebras (\cite{VV}).
Let $a_1, a_2, \cdots , a_d,$ and $q$ be positive integers and we put
$$
Q=(x_1^{a_1}, x_2^{a_2}, \cdots , x_d^{a_d}) \ \ \textup{and} \ \ I=Q:\fkm^q.
$$
Let $\overline{A} = A/Q$, $\overline{\fkm} = \fkm /Q$, and $\overline{I} = I/Q$. Then $$\rmG(\overline{\fkm}) \cong \rmG(\fkm)/({X_1}^{a_1}, {X_2}^{a_2}, \cdots, {X_d}^{a_d}), $$ whence ${\rmG}(\overline{\fkm})$ is a Gorenstein ring. Let $\rho = \max~\{n \in \Bbb Z \mid \overline{\fkm}^n \ne (0)\},$
that is index of nilpotency of the ideal $\overline{\fkm}$, and we have  $$\rho = {\rma}({\rmG} (\overline{\fkm})) = {\rma}({\rmG} (\fkm)) + \sum_{i=1}^da_i,$$ 
where ${\rma}(*)$ denotes the $a$-invariant of the corresponding graded ring (\cite[(3.1.4)]{GW}).

Let $$\ell = \rho + 1 - q.$$ By \cite{Wat} (see \cite[Theorem 1.6]{O} also) we then have the following.

\begin{prop}\label{2.1}
$(0): \overline{\fkm}^i = \overline{\fkm}^{\rho + 1 - i}$ for all $i \in \Bbb{Z}$. In particular $\overline{I} = (0) : \overline{\fkm}^q = \overline{\fkm}^{\ell}$ whence $I = Q + \fkm^{\ell}$.
\end{prop}

The key for our proof of  Theorem \ref{1.2}  is the following.

\begin{lem}\label{2.2}
Suppose that $\ell \ge a_i$ for all $1 \leq i \leq d$. Then 
$$
Q \cap \fkm^{n\ell + m} \subseteq \fkm^m QI^{n-1}
$$ for all $m \ge 0$ and $n \ge 1$. 
\end{lem}

\begin{proof}
We have 
$$
Q \cap \fkm^{n\ell + m} = \sum_{i=1}^d x_i^{a_i} \fkm^{n\ell + m - a_i}
$$
since $x_1, x_2, \cdots, x_d$ is a super regular sequence with respect to $\fkm$ (that is, their initial forms $X_1, X_2, \cdots , X_d$ constitute a regular sequence in $\rmG(\fkm)$). Because $$n\ell + m - a_i = (n-1)\ell + m + (\ell - a_i) \ge (n-1)\ell + m$$  for each  $1 \le i \le d$, we get 
$$\fkm^{n\ell + m - a_i} \subseteq \fkm^{(n-1)\ell + m} = \fkm^m {\cdot}(\fkm^\ell)^{n-1}.$$ Therefore, since  $\fkm^{\ell} \subseteq I$ by Proposition \ref{2.1}, we have   
\begin{eqnarray*}
Q \cap \fkm^{n\ell + m} &=& \sum_{i=1}^d x_i^{a_i} \fkm^{n\ell + m - a_i}\\
&\subseteq& \sum_{i=1}^d x_i^{a_i} \fkm^m(\fkm^\ell)^{n-1}\\
&\subseteq& \fkm^m Q I^{n-1}
\end{eqnarray*}
as is claimed.
\end{proof}

Let us now prove  Theorem \ref{1.2}.

\begin{proof}[Proof of Theorem $\ref{1.2}$]
(2) $\Rightarrow$ (1) This is well-known. See \cite[Section 7, Theorem 2]{NR}.

(3) $\Rightarrow$ (2) By Proposition \ref{2.1}  we get $\fkm^q I = \fkm^q Q + \fkm^{q+\ell}$, whence $\fkm^{q +\ell} \subseteq Q$,  so that $\fkm^{q+\ell} = Q \cap \fkm^{q+\ell}  \subseteq \fkm^{q}Q$ by Lemma \ref{2.2}, because $\ell \ge a_i$ for all $1 \leq i \leq d$. Thus $\fkm^q I = \fkm^q Q$.

(1) $\Rightarrow$ (3) Let $1 \le i \le d$ be an integer. Then $x_i^{\ell} \in \fkm^\ell \subseteq I \subseteq \overline{Q}$.
We look at the integral equation
$$
(x_i^\ell)^n + c_1(x_i^\ell)^{n-1} + \cdots + c_n = 0
$$
with  $n \gg 0$ and $c_j \in Q^j$. Then  $x_i^{n\ell} \in \sum_{j=1}^n Q^jx_i^{(n-j)\ell}$ and so, thanks to the monomial property of the sequence $x_1, x_2, \cdots, x_d$ (cf. \cite[Exercise 17.13, c.]{E}; recall that the sequence $x_1, x_2, \cdots , x_d$ is $A$-regular), we have 
$$x_i^{n\ell} \in Q^jx_i^{(n-j)\ell}$$
for some $1 \leq j \leq n$. Let $L =\{(\alpha_1, \alpha_2, \cdots , \alpha_d) \in {\Bbb Z}^d \mid \alpha_k \ge 0~\textup{for~all}~ 1 \leq k \leq d \}$ and $\Lambda_j = \{\alpha \in L \mid \sum_{k=1}^d \alpha_k = j \}$. Then since  $$Q^j = (\prod_{k=1}^d{x_k}^{a_k\alpha_k} \mid \alpha \in \Lambda_j),$$  thanks to the monomial property of $x_1, x_2, \cdots, x_d$ again, we get 
$$x_i^{n\ell} \in [\prod_{k=1}^{d}{x_k}^{a_k\alpha_k}]{\cdot}x_i^{(n-j)\ell}A$$
for some $\alpha \in \Lambda_j$. Hence $$(n\ell) \mathbf{e}_i = \sum_{k=1}^{d}(a_k\alpha_k)\mathbf{e}_k + [(n-j) \ell]\mathbf{e}_i + \beta$$ with $\beta \in L$, where $\{\mathbf{e}_i\}_{1 \le i \leq d}$ denotes the standard basis of ${\Bbb Z}^d$. Consequently $\alpha_k = \beta_k = 0$ if $k \ne i$ and so $\alpha_i = j$. Hence $0 = a_ij - j \ell + \beta_i$, so that we have $\ell \ge a_i$  as is claimed.

Let us now consider assertions (i) and (ii).
Let $n \geq 1$ be an integer. Then $I^n = QI^{n-1} + \fkm^{n\ell}$ since $I = Q + \fkm^{\ell}$ (Proposition \ref{2.1}), so that $$Q \cap I^n = QI^{n-1} + [Q \cap \fkm^{n\ell}] \subseteq QI^{n-1}$$ because $Q \cap \fkm^{n\ell} \subseteq QI^{n-1}$ by Lemma \ref{2.2}. Therefore $Q \cap I^n = QI^{n-1}$ for all $n \ge 1$, whence $\rmG(I)$ is a Cohen-Macaulay ring  (\cite[Corollary 2.7]{VV}).

We will show that $\rmr_Q(I) = \lceil \frac{q}{\ell} \rceil$. Notice that $$\rmr_Q(I) = \min\{n \ge 0 \mid I^{n+1} \subseteq Q\},$$ because $I^{n+1} = QI^n$ if and only if $I^{n+1} \subseteq Q$. Firstly, suppose that $I^{n+1} \subseteq Q$. We then have $\overline{\fkm}^{(n+1)\ell} = (0)$ (recall that $\overline{I} = \overline{\fkm}^\ell$), whence $(n+1)\ell \ge \rho + 1$. Therefore $$n+1 \ge \frac{\rho +1}{\ell} = \frac{q + \ell}{\ell} = \frac{q}{\ell} +1,$$ because $\ell = \rho +1 - q$, so that we have $n \ge \frac{q}{\ell}$.

If $n \ge \frac{q}{\ell}$, then $(n+1) \ell  \ge (\frac{q}{\ell} +1)\ell = q  + \ell = \rho+1$ and so  $\overline{I}^{n+1} = \overline{\fkm}^{(n+1)\ell} = (0)$, whence $I^{n+1} \subseteq Q$. Thus  $\rmr_Q(I) = \lceil \frac{q}{\ell} \rceil$.

Because $\rmG (I)$ is a Cohen-Macaulay ring, to see the fiber cone $\rmF(I)$ is Cohen-Macaulay it suffices  to show that 
$$
Q \cap \fkm I^n = \fkm QI^{n-1}
$$
for all $n \ge 1$ (see, e.g., \cite{CGPU, CZ}).
By Lemma \ref{2.2} we have
\begin{eqnarray*}
Q \cap \fkm I^n &=& Q \cap [\fkm QI^{n-1} + \fkm^{n\ell +1}]\\
&=& \fkm QI^{n-1} + [Q \cap \fkm^{n\ell + 1}]\\
&\subseteq& \fkm Q I^{n-1},
\end{eqnarray*}
whence $Q \cap \fkm I^n = \fkm Q I^{n-1}$.
\end{proof}

Assume that $\ell \ge a_i$ for all $1 \le i \le d$ and let $Y_i$'s be the initial forms of $x_i^{a_i}$'s with respect to $I$. Then $Y_1, Y_2, \cdots, Y_d$ is a homogeneous system of parameters of $\rmG (I)$, since $Q$ is a reduction of $I$ (Theorem \ref{1.2}). It therefore constitutes a regular sequence in $\rmG (I)$, because $\rmG (I)$ is a Cohen-Macaulay ring by Theorem \ref{1.2} (ii), so that we have a canonical isomorphism 
$$\rmG (\overline{I}) \cong \rmG (I)/(Y_1, Y_2, \cdots, Y_d)$$
of graded $A$-algebras (\cite{VV}). Hence $\rma (\rmG (\overline{I}))= \rma (\rmG (I)) +  d$. Let $r$ be the index of nilpotency of $\overline{I}$, that is $r = \rma (\rmG (\overline{I}))$. Then since  $r = \rmr_Q(I)$ (recall that ${x_1}^{a_1}, {x_2}^{a_2}, \cdots, {x_d}^{a_d}$ is a super regular sequence with respect to $I$) and $\rma (\rmG (I)) = \rma (\rmG (\overline{I})) - d$ (\cite[(3.1.6)]{GW}),  by Theorem \ref{1.2} (i) we have the following.

\begin{lem}\label{a}
Suppose that $\ell \ge a_i$ for all $1 \le i \le d$. Then $\rma (\rmG (I)) = \lceil \frac{q}{\ell} \rceil -d.$
\end{lem}

\begin{cor}\label{CMRees}
Assume that $\ell \ge a_i$ for all $1 \le i \le d$. 
Then $\calR(I)$ is a Cohen-Macaulay ring if and only if $\lceil \frac{q}{\ell} \rceil < d$. When this is the case, $d \ge 2$.
\end{cor}

\begin{proof}
Since $\rmG(I)$ is a Cohen-Macaulay ring by Theorem \ref{1.2} (ii), $\calR(I)$ is a Cohen-Macaulay ring if and only if $a(\rmG(I)) < 0$ (\cite{TI}).
By Lemma \ref{a} the latter condition is equivalent to saying that $\lceil \frac{q}{\ell} \rceil < d$ (cf. \cite[Remark (3.10)]{GSh}). When this is case, $d \ge 2$ because $0 < \lceil \frac{q}{\ell} \rceil$.
\end{proof}

We are now in a position to prove Theorem \ref{1.2'}.

\begin{proof}[Proof of Theorem $\ref{1.2'}$]

(i) Notice that 
$\rmG(I)$ is a Gorenstein ring if and only if so is the graded ring $${\rmG}(\overline{I}) = {\rmG}(I) / (Y_1, Y_2, \cdots, Y_d),$$ where $Y_i$'s stand for the initial forms of $x_i^{a_i}$'s with respect to $I$. 
Let $r$ be the index of nilpotency of $\overline{I}$. Then $r = \rmr_Q(I) = \lceil \frac{q}{\ell} \rceil$, and $\rmG(\overline{I})$ is a Gorenstein ring if and only if the equality $$(0):\overline{I}^i = \overline{I}^{r+1-i}$$ holds true for all $i \in \Bbb{Z}$ (\cite[Theorem 1.6]{O}). Hence if $\rmG(I)$ is a Gorenstein ring, we have $(0) : \overline{I} = \overline{I}^r = \overline{\fkm}^{r\ell}$. On the other hand, since $\overline{I} = \overline{\fkm}^{\ell}$ and $q = \rho + 1 - \ell$, by Proposition \ref{2.1} we get $$(0):\overline{I} = (0):\overline{\fkm}^\ell = \overline{\fkm}^q.$$
Therefore $q = r \ell$, since $\overline{\fkm}^{r\ell} = \overline{\fkm}^{q} \ne (0)$. Thus $\ell \mid q$ and $r =\frac{q}{\ell}$.

Conversely, suppose that $\ell \mid q$. Hence $r = \frac{q}{\ell}$ by Theorem \ref{1.2} (i).
Let $i\in \Bbb{Z}$.
Then since  $\overline{I} = \overline{\fkm}^{\ell}$, we get $\overline{I}^{r+1-i} = \overline{\fkm}^{(r+1-i)\ell}$, while $$(0):\overline{I}^i = (0):\overline{\fkm}^{i\ell} = \overline{\fkm}^{\rho +1 -i\ell}$$ by Proposition \ref{2.1}. We then have $(0):\overline{I}^i = \overline{I}^{r+1-i}$ for all $i \in \Bbb Z$, since $$(r+1-i)\ell = q + \ell - i\ell = \rho +1 - i\ell.$$
Thus $\rmG(\overline{I})$ is a Gorenstein ring, whence so is $\rmG (I)$.

(ii) The Rees algebra $\calR(I)$ of $I$ is a Gorenstein ring if and only if $\rmG(I)$ is a Gorenstein ring and $a(\rmG(I)) = -2$, provided $d \ge 2$ (\cite[Corollary (3.7)]{I}).
Suppose that $\calR(I)$ is a Gorenstein ring. Then $d \ge 2$ by Corollary \ref{CMRees}. Since $a(\rmG(I)) = \rmr_Q(I) -d=-2$, by assertion (i) and  Theorem \ref{1.2} (i) we have $\frac{q}{\ell} = \rmr_Q(I) = d-2$, whence $q= (d-2)\ell$.
Conversely, suppose that $q =(d-2)\ell$. Then $d \ge 3$ since $q \ge 1$. By assertion (i) and Theorem \ref{1.2} (i) $\rmG(I)$ is a Gorenstein ring with $\rmr_Q(I) = \frac{q}{\ell} = d-2$, whence $a(\rmG(I)) = (d-2)-d = -2$, so that $\calR(I)$ is a Gorenstein ring.
\end{proof}

\begin{example}
Suppose that $\rho \ge 5$ is an odd integer, say $\rho = 2\tau + 1$ with $\tau \ge 2$. Let $q = \rho -1$. Then $\ell = \rho + 1 - q= 2$. Hence, choosing $a_i \le 2$ for all $1 \le i \le d$, we have $I = Q + \fkm^2 \subseteq \overline{Q}$ with $\rmr_Q(I) = \tau$ by Theorem \ref{1.2}.  Since $\ell \mid q$, by Theorem \ref{1.2'} (i) the ring ${\rmG}(I)$ is  Gorenstein, despite $Q \not\subseteq \fkm^q$ (compare the condition with those in Theorem \ref{1.4}). The ring $\calR (I)$ is by Theorem \ref{1.2'} (ii) a Gorenstein ring, if $d = \tau + 2$. 
\end{example}


\section{Examples and applications}
In this section we shall discuss some applications of Theorems \ref{1.2} and \ref{1.2'}. Let us begin with the case where $A$ is a regular local ring.

\subsection{The case where $A$ is a regular local ring}

Let $A$ be a regular local ring with $x_1, x_2, \cdots, x_d$ a regular system of parameters. Similarly as in the previous sections, let 
$$
Q=(x_1^{a_1}, x_2^{a_2}, \cdots , x_d^{a_d}) \ \ \textup{and} \ \ I=Q:\fkm^q
$$
with  positive integers $a_1, a_2, \cdots , a_d,$ and $q$. Then ${\rmG}(\fkm)=k[X_1, X_2, \cdots, X_d]$ is the polynomial ring, where $k = A/\fkm$ and $X_i$'s are  the initial forms of $x_i$'s, so that we have $$\rho = \sum_{i=1}^da_i - d \ \ \textup{and} \ \ \ell = \sum_{i=1}^d(a_i-1)  + 1 - q,$$ since $\rma ({\rmG}(\fkm)) = -d$. 
Notice that the condition that $$\ell \geq \max~\{a_i \mid 1 \leq i \leq d\}$$ is equivalent to saying that 
$$\sum_{j\ne i }a_j \geq q + d - 1$$
for all $1 \leq i \leq d$, because $\ell - a_i = \sum_{j \ne i }a_j - (q + d -1)$. When this is the case, $d \ge 2$.

We readily get, thanks to Theorem \ref{1.2},  the following.

\begin{example} The following assertions hold true.
\begin{enumerate}
\item[$(1)$] Let $d = 2$. Then $I \subseteq \overline{Q}$ if and only if $\min \{a_1, a_2 \} \geq q+1$.
\item[$(2)$] Let $d = 3$. Then $I \subseteq \overline{Q}$ if and only if $\min \{a_i + a_j \mid 1 \leq i < j \leq 3\} \geq q + 2$.
\item[$(3)$] Choose integers  $a$ and $q$ so  that $2 \leq a \leq d$ and $(d-1)(a-1) < q \leq d(a-1)$. Let $a_i = a$ for all $1 \leq i \leq d$. Then $I \subsetneq A$ but $I \not\subseteq \overline{Q}$. For example, let $d=3, a=2,$ and $q = 3$.  Then $$(x_1^2, x_2^2, x_3^2) : \fkm^3  = \fkm ~~\not\subseteq ~~\overline{(x_1^2, x_2^2, x_3^2)}.$$
\end{enumerate}
\end{example}

Thanks to Theorem \ref{1.2'}, we are able to produce  quasi-socle ideals $I=Q:\fkm^q$ whose Rees algebras are Gorenstein.

\begin{example}\label{3.2}
The following assertions hold true.
\begin{enumerate}
\item[$(1)$] Let $d = 2$ and assume that $I \subseteq \overline{Q}$. Then $\rmG(I)$ is not a Gorenstein ring.
\item[$(2)$] Suppose that $d \geq 3$ and let $n \geq d-1$ be an integer. Let $a_1 = d-1$, $a_i = n$ for all $2 \leq i \leq d$, and $q = (d-2)n$. Then $\calR (I)$ is a Gorenstein ring. 
\item[$(3)$] Suppose that $d = 5$ and let $a_i = 4$ for all $1 \le i \le 5$.
Let $q = 8$.
Then $I \subseteq \overline{Q}$ and $\rmG(I)$ is a Gorenstein ring with $\rmr_Q(I) = 1$, but $\calR(I)$ is not a Gorenstein ring.
\end{enumerate}
\end{example}

\Proof
(1) Suppose that $\rmG (I)$ is a Gorenstein ring and let $q=r\ell $ with $r = \rmr_Q(I)$.
Then $a_1 + a_2 - 1 = \rho +1 = q+ \ell  = \ell(r+1)$, which implies, because $\ell \ge a_i$ for $i = 1,2$ by Theorem \ref{1.2}, that $r=0$.
This is impossible.

(2) Since $\rho = nd -n -1$, we get $\ell = n$. Hence $q = (d-2)\ell$, so that  by Theorem \ref{1.2'} (ii) $\calR (I)$ is a Gorenstein ring.

(3) We have $\rho = 15$ and $\ell = q = 8$.
Hence $I \subseteq \overline{Q}$ with $\rmr_Q(I) = 1$ by Theorem \ref{1.2}. By Theorem \ref{1.2'} $\rmG(I)$ is a Gorenstein ring, but the ring $\calR(I)$ is not, because $q \ne (d-2) \ell$.
\Qed

\begin{remark}
When $A$ is not a regular local ring, the associated graded ring $\rmG (I)$ of $I$ can be Gorenstein, even though $d = 2$.  See Example \ref{3.8}.
\end{remark}

Since the base ring $A$ is regular, the Cohen-Macaulayness in Rees algebras $\calR (I)$ follows from that of associated graded rings $\rmG (I)$ (\cite{L}). Let us note a brief proof in our context.

\begin{prop}\label{CM}
Suppose that $\ell \ge a_i$ for all $1 \le i \le d$. Then the  Rees algebra $\calR (I)$ is a Cohen-Macaulay ring.
\end{prop}

\begin{proof}
By Corollary \ref{CMRees} we have only to show $\lceil \frac{q}{\ell} \rceil < d$. Let $a_k = \max \{a_i \mid 1 \leq i \leq d\}$. Then because $\ell \geq a_k$, we have 
$$\frac{q}{\ell} + 1 = \frac{\rho + 1}{\ell} \leq \frac{\sum_{j=1}^d(a_j - 1) + 1}{a_k} = \sum_{j \ne k}\frac{a_j -1}{a_k} + 1 < d,$$ whence $\lceil \frac{q}{\ell} \rceil <d$ as is wanted.
\end{proof}

Let $L = \{(\alpha_1, \alpha_2, \cdots, \alpha_d) \in {\Bbb Z}^d \mid \alpha_i \ge 0 ~\textup{for~all} 1 \le i \le d\}$. For each $\alpha = (\alpha_1, \alpha_2, \cdots, \alpha_d) \in L$  we put $x^\alpha = \prod_{i=1}^dx_i^{\alpha_i}$. Let $\fka$ be an ideal in $A$. Then we say that $\fka$ is a monomial ideal, if $\fka$ is generated by monomials in $\{x_i\}_{1 \leq i \leq d}$, that is $\fka = (x^{\alpha} \mid \alpha \in \Lambda)$ for some $\Lambda \subseteq L$. Monomial ideals 
 behave very well as if they were monomial ideals in the polynomial ring $k[x_1, x_2, \cdots, x_d]$ over a field $k$ (see \cite{HS} for details). For instance, the integral closure  $\overline{Q}$  of our monomial ideal $Q$ is also a monomial ideal and  we have the following.

\begin{prop}[\cite{HS}]\label{HS} Let $\Delta = \{\alpha \in L \mid \sum_{i=1}^d\frac{\alpha_i}{a_i} \geq 1\}$. Then $\overline{Q} =(x^{\alpha} \mid \alpha \in \Delta )$. 
\end{prop}

\begin{cor}\label{4.1}
Suppose  that $d \ge 2$ and let $n \geq 2$ be an integer. We put 
$\fkq =(x_1^{n-1}, x_2^n, \cdots, x_d^n)$.
Then $\overline{\fkq } = \fkq + \fkm^n =(x_1^{n-1}) + \fkm^n$ and all the powers $\overline{\fkq}^m$~ $(m \geq 1)$ are integrally closed.
\end{cor}

\Proof
Let $J = \fkq + \fkm^{n}$ and $\fka = (x_1^n, x_2^n, \cdots,  x_d^n)$.
Then $\fka \subseteq \fkq $ and $\fkm^n \subseteq \overline{\fka}$, so that $J \subseteq \overline{\fkq }$. Let $m \geq 1$ be an integer and put $K = (x_1^{m(n-1)}, x_2^{mn}, \cdots, x_d^{mn})$. We will show that $\overline{K} \subseteq J^m$.  Let $\alpha \in L$ and assume that $\frac{\alpha_1}{m(n-1)}+ \sum_{i=2}^d\frac{\alpha_i}{mn} \geq 1$. We want to show that $x^{\alpha} \in J^m$. We may assume  that $\alpha_1 < m(n-1)$. Let $\alpha_1 = (n-1)i + j$ with $i, j \in \Bbb Z$ such that $0 \leq j < (n-1)$. Then $0 \leq i < m$. Since  $\frac{\alpha_1}{m(n-1)}+ \sum_{i=2}^d\frac{\alpha_i}{mn} \geq 1$, we get 
$$n \alpha_1 + (n-1){\cdot}\sum_{i=2}^d\alpha_i \geq mn(n-1), $$  so that 
$$(n-1){\cdot}\sum_{i=2}^d\alpha_i \geq mn(n-1) -n \alpha_1 = n[(n-1)(m-i) -j],$$
 whence $$\sum_{i=2}^d\alpha_i \geq n(m-i) - \frac{nj}{n-1}.$$  Because $\frac{nj}{n-1} = j + \frac{j}{n-1}$ and $0 \leq j < n - 1$, we have $\frac{nj}{n-1} = j + \frac{j}{n-1} < j+1$ and so $$\sum_{i=2}^d\alpha_i \geq n(m-i) -j.$$ Thus  $$x^{\alpha} = x_1^{(n-1)i}{\cdot}x_1^jx_2^{\alpha_2}\cdots x_d^{\alpha_d} \in x_1^{(n-1)i}\fkm^{n(m-i)} \subseteq J^m,$$ whence $\overline{K} \subseteq J^m$ by Proposition \ref{HS}. 

Because $J^m \subseteq {\overline{\fkq}}^m$ and $\fkq^m \subseteq \overline{K}$, we have $J^m \subseteq {\overline{\fkq}}^m \subseteq  \overline{\fkq^m} \subseteq  \overline{K}$, whence $J^m = \overline{\fkq}^m = \overline{\fkq^m} = K$. Letting $m=1$, we get $J = \overline{\fkq}$. This completes the proof of Corollary \ref{4.1}.
\Qed

Thanks to Corollary \ref{4.1}, we  get the following  characterization for quasi-socle ideals $I = Q : \fkm^q$ to be integrally closed.

\begin{prop}\label{4.2}
Suppose that $d \ge 2$ and $a_i \ge 2$ for all $1 \le i \le d$.
Then the following two conditions are equivalent to each other.
\begin{enumerate}
\item[$(1)$] $I=\overline{Q}$.
\item[$(2)$] Either {\rm (a)} $a_i = \ell$ for all $1 \le i \le d$, or {\rm (b)} there exists $1 \le j \le d$ such that $a_i = \ell$ if $i \ne j$ and $a_j = \ell -1$.
\end{enumerate}
When this  is the case, $\overline{I^n} =I^n$ for all $n \geq 1$, whence $\calR (I)$ is a Cohen-Macaulay normal domain.
\end{prop}

\Proof
(1) $\Rightarrow$ (2)
Since $I = \overline{Q}$, we get $q \le \rho$ and $I=Q+\fkm^\ell$ (Proposition \ref{2.1}).
Notice that $$Q \subseteq I =Q:\fkm^q \subsetneq(Q:\fkm^q):\fkm = Q:\fkm^{q+1},$$ because $I \subsetneq A$.
Hence $Q:\fkm^{q+1} \not\subseteq \overline{Q}$.
Consequently $\ell -1 = \rho +1 - (q+1)< a_i$ for some $1 \le i \le d$ by Theorem \ref{1.2}, so that, thanks to Theorem \ref{1.2} again, we have 
$$
\ell = a_i \ge a_j
$$
for all $1 \le j \le d$.
Let $\Delta = \{1 \le j \le d \mid a_j < \ell \}$. We then  have the following.
\begin{claim*}
\begin{enumerate}
\item[$(1)$] $a_j = \ell -1$, if $j \in \Delta$.
\item[$(2)$] $\sharp \Delta \le 1$.
\end{enumerate}
\end{claim*}
\begin{proof}
Let $j \in \Delta$. Then $a_j < \ell = a_i$ whence $j \ne i$ and $\ell \ge 3$.
Let $\alpha = (a_j -1)\mathbf{e}_j + (a_i - a_j)\mathbf{e}_i$.
Then $\alpha \in L$ but, thanks to the monomial property of ideals, $x^\alpha \notin Q+\fkm^\ell = I = \overline{Q}$, because $\sum_{k = 1}^d\alpha_k = a_i -1 = \ell -1$ and $x^\alpha \notin Q$.
Consequently, $\sum_{k=1}^d \frac{\alpha_k}{a_k} < 1$ by Proposition \ref{HS}, so that $1 < \frac{1}{a_j} + \frac{a_j}{a_i}$, because 
$$
\frac{a_j-1}{a_j} + \frac{a_i - a_j}{a_i} < 1.
$$
Let $n=a_i - a_j$.
Then $a_j(a_i - a_j) < a_i$ as $1 < \frac{1}{a_j} + \frac{a_j}{a_i}$, whence $a_j n < a_i = a_j + n$ so that $0 \le (a_j - 1)(n-1) < 1$.
Hence $n=1$ (recall that $a_j \ge 2$) and $a_j = a_i - 1 = \ell -1$.

Assume $\sharp \Delta \ge 2$ and choose $j, k \in \Delta$ so that $j \ne k$.
We put $y = x_j x_k^{\ell-2}$.
We then have $y^{\ell-1} = (x_j^{\ell-1})(x_k^{\ell-1})^{\ell-2} = (x_j^{a_j})(x_k^{a_k})^{\ell-2} \in Q^{\ell-1}$, because $a_j = a_k = \ell -1$ by assertion (1).
Hence $y \in \overline{Q} = Q+ \fkm^\ell$, which is impossible because $y\notin Q$ (recall that $\ell \ge 3$) and $y \notin \fkm^\ell$, thanks to the monomial property of ideals.
Hence $\sharp \Delta \le 1$.
\end{proof}

If $\Delta = \emptyset$, we then have $\ell = a_j$ for all $1 \le j \le d$.
If $\Delta \ne \emptyset$, letting $\Delta = \{j\}$, we get $a_i = \ell$ if $i \ne j$ and $a_j = \ell -1$. This proves the implication $(1) \Rightarrow (2)$.

(2) $\Rightarrow$ (1) Suppose condition (b) is satisfied.
Then $I=Q+\fkm^\ell = (x_j^{\ell -1}) + \fkm^\ell = \overline{Q}$ by Proposition \ref{2.1} and Corollary \ref{4.1}. Suppose condition (a) is satisfied.
Then $I \subseteq \overline{Q}$ by Theorem \ref{1.2} and $I=Q+\fkm^\ell = \fkm^\ell$ by Proposition \ref{2.1}, whence $I=\overline{Q}$. In each case all the powers of $I$ are integrally closed (see Corollary \ref{4.1} for case (b)),
whence the last assertion follows from Proposition \ref{CM}. 
\Qed

\begin{example}
Suppose that $d \ge 3$ and let $n \ge d - 1$ be an  integer. We look at the ideal $$Q = (x_1^{d-1}, x_2^n, x_3^n, \cdots, x_d^n)$$ and let $q = n(d-2)$. Then $\ell = n$, as $\rho = nd - (n + 1)$, whence $I \subseteq \overline{Q}$ and  $I = Q + \fkm^n = (x_1^{d-1}) + \fkm^n$. The ring $\calR (I)$ is by Theorem \ref{1.2'} (ii) a Gorenstein ring, since $q = (d-2)\ell$. If $n = d$, then $I = (x_1^{d-1}) + \fkm^d$ and $\overline{I^m} = I^m$ for all $m \ge 1$ by Corollary \ref{4.1}, so that  $\calR (I)$ is a Gorenstein normal ring.
\end{example}

\subsection{The case where $A = R_M$}
Our setting naturally contains the case where $A = R_M$ is the localization of the homogeneous Gorenstein ring $R= k[R_1]$ over an infinite field $k = R_0$ at the irrelevant maximal ideal $M = R_+$. Let us note one example.

\begin{example}\label{3.8}
Let $S = k[X, Y, Z]$ be the polynomial ring over an infinite field $k$ and let $R = S/fS$, where $0 \ne f \in S$ is a form with degree $n \ge 2$. Then $R$ is a homogeneous Gorenstein ring with $\operatorname{dim} R = 2$. Let $x_1, x_2$ be a linear system of parameters in $R$ and let $M = R_+$. We look at the local ring $A = R_M$. Let $a_1 = 2$, $a_2 = n$, and $q = n$. Let $Q = (x_1^2, x_2^{n})A$ and $I = Q : \fkm^q$, where $\fkm = MA$. Then $$\rho = \rma (R) + (a_1 + a_2) = 2n-1.$$ Hence $\ell = q = n$, so that $I \subseteq \overline{Q}$, $I = Q + \fkm^n = (x_1^2) + \fkm^n$, and ${\rmG}(I)$ is a Gorenstein ring with $\rmr_Q(I) = 1$ (Theorems \ref{1.2} and \ref{1.2'}). We have $Q \not\subseteq \fkm^q$, if $n \ge 3$. 
\end{example}

\subsection{The case where $A = k[[t^a, t^b]]$}
Let $1 < a <b$ be integers with $\GCD(a,b) = 1$. We look at the ring $A = k[[t^a, t^b]] \subseteq k[[t]]$, where $k[[t]]$ denotes the formal powers series ring over a field $k$. We put $x = t^a$ and $y = t^b$.
Then $A$ is a one-dimensional Gorenstein local ring and $\fkm = (x, y)$. Because $A \cong k[[X,Y]]/(X^b-Y^a)$ where $k[[X,Y]]$ denotes the formal powers series ring over the field $k$, we get $${\rmG} (\fkm) \cong k[X,Y]/(Y^a).$$ 

Let $n, q \ge 1$ be integers, and put $Q = (x^n)$ and $I = Q: \fkm^q$. Then because  $a(\rmG(\fkm)) = a - 2$, we have $\rho = a+n-2$ and $\ell = (a+n) - (q+1)$. Consequently $I \subseteq \overline{Q}$ if and only if $q < a$ (Theorem \ref{1.2}), whence the condition that $I \subseteq \overline{Q}$ is independent of the choice of the  integer $n \ge 1$. When this is the case, by Theorems \ref{1.2} and \ref{1.2'} we have the following.

\begin{prop} The following assertions hold true.
\item[(1)] $\rmr_Q(I) = \lceil \frac{q}{ (a+n) - (q+1)} \rceil$.
\item[(2)] The graded rings $\rmG(I)$ and $\rmF(I)$ are Cohen-Macaulay rings.
\item[(3)] The ring $\rmG(I)$ is a Gorenstein ring if and only if $ (a+n) - (q+1)$ divides $q$.
\end{prop}

\noindent
Hence, if $q = a-1$, we then have, for each integer $n \ge 1$  such that $n \mid q$, that  $\rmG (I)$ is a Gorenstein ring. 

\end{document}